\documentclass{irmadegm}
\usepackage[psamsfonts]{amssymb}
\usepackage{amsmath,cmmib57}
\usepackage{amsthm}
\usepackage{latexsym}

\makeatletter
\def\thmhead@plain#1#2#3{%
  \thmname{#1}\thmnumber{\@ifnotempty{#1}{ }#2}%
  \thmnote{ \the\thm@notefont(#3)}}
\let\thmhead\thmhead@plain
\def\swappedhead#1#2#3{%
  \thmnumber{#2}\thmname{\@ifnotempty{#2}{. }#1}%
  \thmnote{ \the\thm@notefont(#3)}}
\thm@notefont{\rm} \makeatother

\theoremstyle{definition} 

 \newtheorem{definition}{Definition}[section]
 
 \newtheorem{example}[definition]{Example}

\theoremstyle{plain}      

 \newtheorem{proposition}[definition]{Proposition}
 \newtheorem{theorem}[definition]{Theorem}
 
 \newtheorem{lemma}[definition]{Lemma}

\newcommand{\De}{\Delta}
\newcommand{\ep}{\varepsilon}
\newcommand{\vfi}{\varphi}
\newcommand{\la}{\lambda}
\newcommand{\tr}{\mathrm{tr}}
\newcommand{\ot}{\otimes}
\newcommand{\ebew}{\mbox{} \hfill{$\scriptscriptstyle{\blacksquare}$}}

\begin{document}

\title{Multiplier Hopf~{$^*$-algebras} with positive integrals: A laboratory for locally compact quantum groups}

\author{Alfons Van Daele}

\address{Department of Mathematics, K.U.\ Leuven \\
Celestijnenlaan~200B, B-3001 Heverlee, Belgium\\
email:\,\texttt{Alfons.VanDaele@wis.kuleuven.ac.be}}

\markboth{Alfons Van Daele} {Multiplier Hopf~{$^*$-algebras} with positive integrals}

\maketitle

\begin{abstract} Any multiplier Hopf~{$^*$-algebra} with positive
integrals gives rise to a locally compact quantum group (in the sense
of Kustermans and Vaes). As a special case of such a situation, we have
the compact quantum groups (in the sense of Woronowicz) and the discrete
quantum groups (as introduced by Effros and Ruan). In fact, the class
of locally compact quantum groups arising from such multiplier
Hopf~{$^*$-algebras} is self-dual.

The most important features of these objects are (1) that they are of a
purely algebraic nature and (2) that they have already a great
complexity, very similar to the general locally compact quantum groups.
This means that they can serve as a good model for the general objects,
at least from the purely algebraic point of view.  They can therefore
be used to study various aspects of the general case, without going
into the more difficult technical aspects, due to the complicated analytic
structure of a general locally compact quantum group.

In this paper, we will first recall the notion of a multiplier
Hopf~{$^*$-algebra} with positive integrals. Then we will illustrate how
these {\it algebraic quantum groups} can be used to gain a deeper
understanding of the general theory.
An important tool will be the Fourier transform. We will also
concentrate on certain actions and how they behave with respect to this
Fourier transform. On the one hand, we will study this in a purely
algebraic context while on the other hand, we will also pass to the
Hilbert space framework.
\end{abstract}

\section{Introduction}

A locally compact quantum group is a pair $(A,\De)$ of a {$C^*$-algebra}
$A$ and a comultiplication $\De$ on $A$, satisfying certain properties.
If $A$ is an abelian {$C^*$-algebra}, it then has the form $C_0(G)$, the
{$C^*$-algebra} of all continuous complex functions, tending to~$0$ at
infinity on a locally compact group $G$ and the comultiplication $\De$
is given by the formula $(\De(f))(p,q)=f(pq)$ where $f\in C_0(G)$ and
$pq$ is the product in $G$ of the elements $p,q$. Observe that in this
case, $\De(f)$ is a bounded continuous complex function on $G\times G$
that in general, will not belong to $C_0(G\times G)$. Indeed, also in
the general case, the comultiplication $\De$ is a $^*$-homomorphism on
$A$ with values in $M(A\ot A)$, the multiplier algebra of the spatial
{$C^*$-tensor} product $A\ot A$ of $A$ with itself.

Any locally compact group $G$ carries a left and a right Haar measure.
This is also true for a locally compact quantum group. In this case,
these are (nice) invariant weights on the {$C^*$-algebra}. An important
fact is that, in the quantum case, the existence of these weights is
part of the axioms, whereas in the classical theory, it is possible to
prove the existence of the Haar measures. Such existence theorems exist
only in special cases for locally compact quantum groups however (and
it seems that a general existence theorem is still out of sight).

The structure of a locally compact quantum group is very rich, but also
technically difficult to work with. It requires not only the standard
results on operator algebras (like the Tomita-Takesaki theory), but
also fundamental skills with weights on {$C^*$-algebras}, unbounded
operators on Hilbert spaces, $\ldots$ And all of this comes on top of a
highly non-trivial algebraic structure, involving a lot of objects.
This makes it rather hard to work with locally compact quantum groups.
Moreover, the present non-trivial examples are very complicated
(although something is changing here, thanks to work done by Vaynerman
and Vaes, see e.g.\ \cite{[V-V1]} and also \cite{[V-V2]}).
All of this makes it difficult, and perhaps
not very attractive, to try to learn the theory and start working in
it. Nevertheless, all the people familiar with the theory know that the
structure is very rich and that this is a nice piece of mathematics.

Fortunately, there are the `algebraic quantum groups'. These are
multiplier Hopf~{($^*$-)algebras} with (positive) integrals (see
section~\ref{sec2} where we start by recalling this notion). As we
mentioned already in the abstract, any multiplier Hopf~{$^*$-algebra}
with positive integrals
gives rise (in a straightforward and easy way) to a locally compact
quantum group. However, not all locally compact quantum groups are of
this form. The compact and discrete quantum groups belong to this class
and some combinations of those two (like the Drinfel'd double of a
compact quantum group). The class is also self-dual. Among the locally
compact  groups it seems to be possible to characterize those coming
from a multiplier Hopf algebra (\cite{[L]}). Such a result is not yet
available for locally compact quantum groups.

Algebraic quantum groups are of a purely algebraic
nature and it is possible to work with them without going into deep
analysis. Nevertheless, the structure is very rich and from the
algebraic point of view, contains all features of the general locally
compact quantum groups. All of the relevant data are present and
essentially no extra relations are imposed by the restriction to these
algebraic quantum groups. For completeness, we have to mention however
that we still don't know of examples of algebraic quantum groups where
the scaling group is not leaving the integrals invariant -- this is
still open. On the other hand, in the non $^*$-case, such examples are
known (and are in fact not so complicated), see e.g.\ \cite{[VD3]}
or \cite{[VD-Z]}.

It seems fair to say that the development of the general theory of
locally compact quantum groups (by Kustermans and Vaes, see
\cite{[K-V1]} and \cite{[K-V2]}) has been
possible, among other reasons, because of the work done before by
Kustermans and myself on algebraic quantum groups (see
\cite{[K-VD1]}). And, as we indicated
already above, it is a common practice to verify general results about
locally compact quantum groups first for algebraic quantum groups
(where only the algebraic aspects have to be considered). Extending
these results to the general case later is usually fairly complicated,
but there is always a good chance that it can be done. Indeed,
algebraic quantum groups are {\it a good model} for general locally compact
quantum groups.

And there is more. Not only to obtain new results, but also for
understanding the old ones, it is important to get first some
familiarity with the framework of algebraic quantum groups (as
probably, the authors themselves have done before obtaining their
general result; of course, not publishing this intermediate step).

This is precisely what this note is all about: After recalling some of
the basics of multiplier Hopf~{$^*$-algebras} with positive integrals (in
section~\ref{sec2}), we illustrate the above strategy in the two
following sections. In section~\ref{sec3}, we take a certain point of
view, starting from a dual pair
of multiplier Hopf~{$^*$-algebras}. The $^*$-structure however does not
play an essential role here. On the other hand, in section~\ref{sec4}, we pass
to the Hilbert space level and there the $^*$-operation and positivity of
the integrals becomes essential.

The key to our approach here is the Fourier transform. In the general
theory, the Fourier transform is not very explicit. The main reason
is that the Hilbert spaces, $L^2(G)$ and $L^2(\hat G)$ in the
classical case of an abelian locally compact group, are identified
through this Fourier transform in the general quantum case. This common
practice has clear advantages, but it also makes some features less
transparent.

In the present note, very few proofs are given. In the first part of
this paper we recall some of the basic notions and known results. Details can
be found elsewhere and references will be given. On the other hand,
many other results that we present later, are not yet found (in this
form) in other papers and it is our intention to publish details
together with J.\ Kustermans in \cite{[K-VD2]}. However, we must say that
essentially most of the results are, in some form, already
present in one of the papers \cite{[K-V2]}, \cite{[K-V3]}
and \cite{[K-VD1]}. The main
difference is the explicit use of the Fourier transform. Recall that after
all, this paper is meant to serve mainly as a tool for learning and
understanding the subject.

For the standard notions and results on Hopf algebras, we refer to the
basic works of Abe \cite{[A]} and Sweedler \cite{[Sw]}. For some
information about dual pairs of Hopf algebras, we refer to
\cite{[VD1]}. For the theory of multiplier Hopf algebras, the reference
is \cite{[VD2]} while algebraic quantum groups (multiplier Hopf algebras
with integrals) are studied in \cite{[VD3]}. Dual pairs of multiplier
Hopf algbras are treated in \cite{[D-VD]} and actions of multiplier
Hopf algebras in \cite{[D-VD-Z]}. A survey on the theory of multiplier
Hopf algebras is given in \cite{[VD-Z]}.

Then, as part of this work also takes place in Hilbert spaces, we need
to give some references about operator algebras also. Much of this can
be found in \cite{[K-R]} but also a good reference is \cite{[S-Z]}.
For the theory of weights and the Tomita-Takesaki theory, we refer
to \cite{[St]}. The theory of Kac algebras is to be found in \cite{[E-S]}.

Finally, we would like to say something about conventions. The algebras
we deal with are algebras over the complex numbers and may or may not
have an identity. If there is no identity however, the product is
assumed to be non-degenerate (as a bilinear form). We are mainly
interested in {$^*$-algebras}. These are algebras with an involution
$a\mapsto a^*$ satisfying the usual properties. Essentially, these
{$^*$-algebra} structures are always of a certain type because we assume
that there is a faithful positive linear functional. For the
comultiplications, we use the symbol $\De$. This comes from Hopf
algebra theory. However, this choice is not completely obvious here as
the same symbol is also commonly used for the modular function of a
non-unimodular locally compact group and (related) for the modular
operator in the Tomita-Takesaki theory. We will use other symbols for
these objects.

In fact, the difficulty arises from the fact that this
material is relating two completely different fields in mathematics.
The first one is the theory of Hopf algebras and the second one is the
theory of operator algebras. Different customs are usual in these two
areas. Since we are mainly interested in the theory of locally compact
quantum groups, that is formulated in the operator algebra framework,
we will follow what is common there. We will use however the Sweedler
notation as it is justified to do so and of course it makes many
formulas and arguments much more transparent. Indeed, we would like
this paper also to be readable for the Hopf algebra people and we
hope that our third section (where we do not emphasize on the
involutive structure) will serve as a bridge between the two areas.\\

\textbf{Acknowledgements.} First, I would like to thank my colleagues (and friends)
at the Institute of Mathematics in Oslo,
where part of this work was done, for their hospitality during my visit
in November 2001. Secondly, I am grateful to the organizers of the
meeting in Strasbourg, in particular L.~Vaynerman, for giving me the
opportunity to talk about my work. I also like to thank my coworkers
{J.~Kustermans} and {S.~Vaes} for many fruitful discussions on this
subject. Finally, I like to thank A.~Jacobs for some \LaTeX-help.

\section{Algebraic quantum groups}\label{sec2}

We will first briefly recall the notion of a multiplier Hopf~{$^*$-algebra}.
For details, we refer to \cite{[VD2]}, see also \cite{[VD-Z]}.

\begin{definition}\label{def2.1}
A {\it multiplier Hopf~{$^*$-algebra}} is a pair $(A,\Delta)$ of a {$^*$-algebra}
$A$ (with a non-degenerate product) and a comultiplication $\Delta$ on
$A$ such that the linear maps $T_1$ and $T_2$ defined on $A\ot A$ by
\begin{align*}
  T_1(a \ot a')&=\Delta(a)(1\ot a')\\
  T_2(a \ot a')&=(a\ot 1)\Delta(a')
\end{align*}
are one-to-one and have range equal to $A\ot A$.
\end{definition}

We have to give some more explanation.

The {$^*$-algebra} may or may not
have an identity. However, the product, as a bilinear map, must be
{\it non-degenerate}. This is automatic if an identity~exists. For an algebra
$A$ with a non-degenerate product, one can define the so-called {\it
multiplier algebra}. It contains $A$ as an essential ideal and it has
an identity. In fact, it is the largest algebra with these properties.
Because $A$ is a {$^*$-algebra}, the multiplier algebra $M(A)$ is also a
{$^*$-algebra}. The tensor product $A\ot A$ is again a {$^*$-algebra} with
a non-degenerate product and also the multiplier algebra $M(A\ot A)$
can be constructed. Elements of the form $1\ot a$ and $a\ot 1$ exist in
$M(A\ot A)$ for all $a\in A$.

A {\it comultiplication} on $A$ is a $^*$-homomorphism $\De:A\to M(A\ot A)$
which is non-degenerate and coassociative. To be non-degenerate here
means that $\De(A)(A\ot A)=A\ot A$. This property is automatic when $A$
has an identity~$1$ and when $\De$ is unital, i.e.\ $\De(1)=1\ot 1$.
Because of the non-degeneracy of~$\De$, it is possible to extend the
obvious maps $\De\ot\iota$ and $\iota\ot\De$ (where $\iota$ is the
identity map) on $A\ot A$ to maps from
$M(A\ot A)$ to $M(A\ot A\ot A)$. This is why coassociativity makes
sense in the form $(\De\ot\iota)\De=(\iota\ot\De)\De$. Finally, the
linear maps $T_1$ and $T_2$, as defined in the definition, will be maps
from $A\ot A$ to $M(A\ot A)$. The requirement is that they are
injective, have range in $A\ot A$ and that all of $A\ot A$ is in the
range of these maps.

The following is the motivating example for this notion.

\begin{example}\label{ex2.5}
Let $G$ be a group and let $A$ be the algebra $K(G)$ of complex
functions with finite support in $G$. Then $A\ot A$ is identified with
$K(G\times G)$ while $M(A\ot A)$ is the algebra of all complex
functions on $G\times G$. The map~$\De$, defined by $\De(f)(p,q)=f(pq)$
whenever $p,q\in G$ and $f\in K(G)$, will be a comultiplication on $A$.
Coassociativity is a consequence of the associativity of the group
multiplication in $G$.
\end{example}

Here is the relation with the notion of a Hopf~{$^*$-algebra}
(see \cite{[VD2]}).

\begin{proposition}
If $(A,\De)$ is a Hopf~{$^*$-algebra}, then
it is a multiplier Hopf~{$^*$-algebra}. Conversely, if $(A,\De)$ is a
multiplier Hopf~{$^*$-algebra} and if $A$ has an identity, then it is a
Hopf~{$^*$-algebra}.
\end{proposition}

 \textbf{\emph{Proof.}} (sketch)\ \
i) If $(A,\De)$ is a Hopf~{$^*$-algebra}, the inverses of the maps $T_1$
and $T_2$ in Definition~\ref{def2.1} are given in terms of the antipode $S$:
\begin{align*}
  T_1^{-1}(a\ot a')&=(\iota\ot S)(\De(a))(1\ot a') \\
  T_2^{-1}(a\ot a')&=(a\ot 1)(S\ot \iota)(\De(a')).
\end{align*}

ii) On the other hand, if $(A,\De)$ is any multiplier Hopf~{$^*$-algebra},
the above formulas can be used to construct an antipode
(and a counit). The counit is a $^*$-homomorphism $\ep:A\to
\mathbb{C}$ such that $(\ep\ot\iota)\De(a)=a$ and $(\iota\ot \ep)\De(a)=a$
for all $a\in A$. The antipode is a anti-homomorphism $S:A\to A$
satisfying $S(S(a)^*)^*=a$ and
\begin{align*}
  m(S\ot\iota)\De(a)&=\ep(a)1\\
  m(\iota\ot S)\De(a)&=\ep(a)1
\end{align*}
for all $a\in A$ (where $m$ is multiplication as a linear map from
$A\ot A$ to $A$). These formulas are given a meaning in $M(A)$.

iii) So, if $(A,\De)$ is a multiplier Hopf~{$^*$-algebra} with an
identity, then it is automatically a Hopf~{$^*$-algebra}.
\ebew\\

It is obvious that the antipode $S$ and the counit $\ep$ in the case of
the Example~\ref{ex2.5} are given by $S(f)(p)=f(p^{-1})$ and $\ep(f)=f(e)$
where $e$ is the identity in the group and where $p^{-1}$ is the
inverse of $p$.\\

Next we recall the notion of an integral on a multiplier Hopf~{$^*$-algebra}
(see \cite{[VD3]}).

\begin{definition}
A linear functional~$\vfi$ on~$A$ satisfying
{$(\iota\ot\vfi)\De(a)=\vfi(a)1$} for all $a\in A$ is called left invariant.
A linear
functional~$\psi$ on~$A$ is called right invariant if
$(\psi\ot\iota)\De(a)=\psi(a)1$ for all $a\in A$. A non-zero left
invariant functional is called a {\it left integral} while a non-zero
right invariant functional is called a {\it right integral}.
\end{definition}

Observe that the above formulas, expressing invariance, again must be
considered in $M(A)$. Left invariance of $\vfi$ should e.g.\ be written
in the form $(\iota\ot\vfi)((a'\ot 1)\De(a))=\vfi(a)a'$ for
all~$a,a'\in A$.

We have the following results on integrals on a multiplier
Hopf~{$^*$-algebra} (see~\cite{[VD3]}).

\begin{proposition}\label{prop2.5}
Let $(A,\De)$ be a multiplier Hopf~{$^*$-algebra} and assume that a left
integral $\vfi$ exists. Then we have:
\begin{enumerate}\setlength{\itemsep}{-\itemsep}
\item[i)] There is also a right integral $\psi$.
\item[ii)] The left and right integrals are unique, up to a scalar.
\item[iii)] The integrals are faithful.
\item[iv)] There exists an invertible multiplier $\delta$ in $M(A)$
such that $(\vfi\ot\iota)\De(a)=\vfi(a)\delta$ and
$(\iota\ot\psi)\De(a)=\psi(a)\delta^{-1}$.
\item[v)] There exists automorphisms $\sigma$ and $\sigma'$ of $A$
such that $\vfi(aa')=\vfi(a'\sigma(a))$ and
$\psi(aa')=\psi(a'\sigma'(a)))$ for all $a,a'\in A$.
\item[vi)] There exists a scalar $\nu\in \mathbb{C}$ such that
$\vfi(S^2(a))=\nu\vfi(a)$ for all $a\in A$.
\end{enumerate}
\end{proposition}

In fact, this result is true for any regular multiplier Hopf
algebra (again see~\cite{[VD3]}).

In the case of a multiplier Hopf~{$^*$-algebra}, it is natural to assume
that the left integral $\vfi$ is positive, i.e.\ that $\vfi(a^*a)\geq
0$ for all $a\in A$. In that case, we have some more consequences for
the above data. It can be shown e.g.\ that automatically a positive
right integral exists. This is not obvious however and the only
argument available now is given in \cite{[K-VD1]}. Also, we have
that~$|\nu|=1$. Furthermore, some very nice analytic properties can be
proven about the multiplier $\delta$ and about the
automorphisms $\sigma$ and $\sigma'$ (see \cite{[K2]}).

There are many extra properties and relations among the different data
that appear in Proposition~\ref{prop2.5}. We refer to \cite{[K-VD2]}
for a collection of them.

From now on, we will assume that $(A,\De)$ is a multiplier Hopf~{$^*$-algebra}
with positive integrals. We will also use the term {\it
{($^*$-)algebraic} quantum group}. Observe that the adjective `{\it
algebraic}'
is not referring to the concept of an algebraic group, but rather to the
purely algebraic framework that one can use in the study of this type
of locally compact quantum groups.

We will use $\vfi$ to denote a
positive left integral and we use $\psi$ for a positive right integral.
It is possible to give a standard relative normalization of~$\vfi$
and~$\psi$ (see e.g. \cite{[K-VD2]}), but we will not need it.\\

We now turn our attention to {\it the dual}. We have the following result
(again see \cite{[VD3]}).

\begin{theorem}
Let $(A,\De)$ be a multiplier Hopf~{$^*$-algebra} with positive
integrals. Let $\hat A$ be the subspace of the dual space  of $A$,
consisting of linear functionals on $A$ of the form $\vfi(\,\cdot\,a)$
where $\vfi$ is a left integral and $a\in A$. Then $\hat A$ can be made
into a multiplier Hopf~{$^*$-algebra} and again it has positive
integrals.
\end{theorem}

First observe that the elements of $\hat A$ are also the ones of the
form $\vfi(a\,\cdot\,)$, $\psi(\,\cdot\,a)$ or $\psi(a\,\cdot\,)$ where
always $a$ runs through $A$.  The product in $\hat A$ is
obtained by dualizing the coproduct and the involution is given by
{$\omega^*(a)=\omega(S(a)^*)^-$} where~$S$ is the antipode. The coproduct
$\hat \De$ on $\hat A$ is dual to the product in $A$. A right integral
$\hat \psi$ on $A$ is obtained by letting $\hat \psi(\omega)=\ep(a)$
when $\omega=\vfi(\,\cdot\,a)$. With this definition, we get
$\hat\psi(\omega^*\omega)=\vfi(a^*a)$ when as before $\omega=
\vfi(\,\cdot\,a)$. So indeed, $\hat \psi$ is again positive.
Applying the procedure once more takes us back to the original
multiplier Hopf~{$^*$-algebra}. For details here see~\cite{[VD3]}.

There are also many formulas relating the data for $(A,\De)$,
namely $\ep$, $S$, $\vfi$, $\psi$, $\sigma$, $\sigma'$, $\nu$ with
the corresponding data for the dual algebra $(\hat A,\hat\De)$.
See e.g.~\cite{[K2]} and also \cite{[K-VD2]}.\\

In the {\it remaining of this paper}, we will let $B$ be the dual
algebra $\hat A$, but we will consider it with the {\it opposite
comultiplication} $\hat\De^{\text{op}}$. We will systematically
use letters $a, a', \ldots$ (and sometimes $x, x'$) to denote
elements in $A$ and letters $b, b', \ldots$ (or $y, y'$) for
elements in $B$.  We will use $\langle a, b\rangle$ for the
evaluation of the element $b$ in the element $a$. We will (in
general) also drop the symbol $\,\hat{}\,$ on the objects related
with $B$. In other words, we will be working with a (modified)
dual pair $(A,B)$ of multiplier Hopf~{$^*$-algebras} with positive
integrals (as introduced in \cite{[D-VD]}). We call it `{\it modified}'
because of the fact that the coproduct in~$B$ has been reversed
(contrary to the original definition in \cite{[D-VD]}).  This has
e.g.\ as a consequence that $S=\hat S^{-1}$ on $B$ and so
\[\langle S(a),b \rangle= \langle a, S^{-1}(b) \rangle\]
for all $a\in A$ and $b\in B$. Also $\vfi=\hat \psi$ and $\psi=\hat\vfi$
on $B$.

The reason for this modification is to get the theory here in
accordance with the general theory of locally compact quantum groups.

A last remark for this section: we will use the Sweedler notation. It is
justified as long as we have the right {\it coverings}: e.g.\
$\De(a)(1\ot a')$ is written as $a_{(1)}\ot a_{(2)}a'$ and we say that
$a_{(2)}$ is covered by $a'$. See \cite{[D-VD]} for a detailed discussion
on the use of the Sweedler notation for multiplier Hopf algebras. Observe
that sometimes, the covering is through the pairing. The reason is that
given an element $b$ in $B$ e.g.\ there exists an element $e$ in $A$
such that $\langle a, b \rangle= \langle ea, b \rangle$ for all $a \in A$
(cf.\ \cite{[D-VD-Z]}) so that the
element $b$ will, in a sense, cover the element $a$, through the
pairing, by means of this element $e$.

\section{Actions and the Fourier transform}\label{sec3}

As in the previous section, also here $(A,B)$ will be a dual pair of
multiplier Hopf~{$^*$-algebras} with
positive integrals. Again we will assume that $B=\hat A$, and that $B$
is endowed with the
opposite coproduct $\De=\hat\De^{\text{op}}$.

We will now define an action of $A$ and an action of $B$ on $A$.

\begin{definition}\label{def3.1}
We define linear maps $\pi(a)$ for $a\in A$
and $\la(b)$ for $b\in B$ on the space $A$ by
\begin{align*} \pi(a)x&=ax\\
      \la(b)x&=\langle S^{-1}(x_{(1)}),b\rangle x_{(2)}
\end{align*}
where $x\in A$.
\end{definition}
Observe that $x_{(1)}$ is covered, through the pairing, by $b$.

One can verify that these are indeed actions of $A$ and of $B$.

Remark that we have used $x$ for an element in $A$. We will do this
systematically for
elements in $A$ when we treat $A$ as the space on which the algebras
act. We will use the letter $y$ for elements in $B$ when we have
actions on the space $B$.

Throughout this paper, the reader should have in mind that
$\pi(\,\cdot\,)$ will be used for `{\it multiplication}' operators and
that $\lambda(\,\cdot\,)$ will be used for `{\it convolution}'
operators.

We have the following commutation rules.

\begin{lemma}\label{lem3.2}
For $a\in A$ and $b\in B$ we have
\[\pi(a)\la(b)=\langle a_{(1)} , b_{(1)}\rangle\, \la(b_{(2)})\pi(a_{(2)}).\]
\end{lemma}

These relations are called the {\it Heisenberg commutation relations}.
One can show that the linear span of the operators $\la(b)\pi(a)$ is
precisely the linear span of the rank one maps on $A$ of the form
$x\mapsto \langle x,b'\rangle a'$ where $a'\in A$ and $b'\in B$.

We have already considered the bijection $a\mapsto \vfi(\,\cdot\, a)$
from $A$ to $B$. This is one possible {\it Fourier transform}. We will
also consider another one:

\begin{definition}\label{def3.3}
We define two maps $F_1,F_2:A\to B$ by
\begin{align*} F_1(a)&=\vfi(\,\cdot\,a)\\
         F_2(a)&=\psi(S(\,\cdot\,)a)
\end{align*}
\end{definition}

Observe that these maps depend on the normalization of $\vfi$ and
$\psi$. They are related by a simple formula: $F_2$ is a scalar
multiple of $F_1\circ \sigma\circ S^{-1}$. The two maps can also be
regarded in their dual forms. Then, we get the following:

\begin{proposition}\label{prop3.4}
For all $a\in A$ and $b\in B$ we have
\begin{align*}
i) \qquad b&=\vfi(\,\cdot\,a) & \quad \Leftrightarrow \quad
a&=\vfi(S^{-1}(\,\cdot\,)b)\\
ii)\qquad b&=\psi(S(\,\cdot\,)a) & \quad \Leftrightarrow \quad
a&=\psi(\,\cdot\,b)
\end{align*}
\end{proposition}

In this case, we need the relative normalization of the left
integrals on~$A$ and~$B$ together with the relative normalization
of the right integrals on~$A$ and~$B$ w.r.t.\ each other (cf.\
\cite{[VD2]} and also \cite{[K-VD2]}).

Later in this section, we will give another form of these formulas.

As we mentioned already, these maps play the role of the Fourier
transforms. Observe that $\vfi(b^*b)=\vfi(a^*a)$ when $b=F_1(a)
=\vfi(\,\cdot\,a)$ (see further).
This is a result that we have mentioned before
already, taken into account that $\vfi=\hat \psi$ on $B$. For the other
transform, we have $\psi(b^*b)=\psi(a^*a)$ when
$b=F_2(a)=\psi(S(\,\cdot\,)a)$. As we have chosen to work basically with the
left integrals, we will be mainly using the Fourier transform $F_1$ and
we will drop the index and simply write $F$ for $F_1$

What happens with the actions under this Fourier transform? The Fourier
transform does convert the multiplication operators to convolution
operators and the other way around. More precisely, we get
the following:

\begin{proposition}\label{prop3.5}
For all $a,x\in A$ and $b\in B$ we have
\begin{align*} F(\pi(a)x)&=\la(a)F(x)\\
         F(\la(b)x)&=\pi(b)F(x)
\end{align*}
where $\la(a)$ and $\pi(b)$ are defined on $B$ by
\begin{align*} \la(a)y&=\langle a,y_{(1)}\rangle y_{(2)} \\
         \pi(b)y&=by
\end{align*}
with $y\in B$.
\end{proposition}

One can verify e.g.\ that, for $a,a'\in A$ and $y\in B$,
\begin{align*} \la(aa')y &= \langle aa', y_{(1)}\rangle y_{(2)} \\
 &=\langle a,y_{(2)}\rangle  \langle a', y_{(1)}\rangle y_{(3)} \\
 &=\la(a) \langle a',y_{(1)}\rangle y_{(2)} \\
 &=\la(a)\la(a')y.
\end{align*}
Observe the difference with the action $\la$ of $B$ on $A$ where the
antipode was involved. This is not so here, for the action $\la$ of $A$
on $B$, a fact that is related with taking the opposite
comultiplication on the dual (as we can notice in the above argument).

We have a similar set of formulas when we look at $F_2$ in stead of
$F_1$.\\

Now, we relate all of this with the so-called {\it left regular
representation}.

\begin{definition}
Consider the map $V$ from $A\ot A$ to itself defined by
\[V(x\ot x') =\De(x')(x\ot 1).\]
\end{definition}

This map is bijective and the inverse is given by
\[V^{-1}(x\ot x')= S^{-1}(x'_{(1)}) x \ot x'_{(2)}\]
(using the Sweedler notation). We will denote $V^{-1}$ by $W$. Then we
get the following:

\begin{lemma}\label{lem3.7}
The map $W$ verifies the Pentagon equation
\[ W_{12} W_{13}W_{23} = W_{23 }W_{12}\]
(where we use the {\it leg-numbering notation}).
\end{lemma}

This equation is to be considered, in the first place, as an equality
of linear maps from $A\ot A\ot A$ to itself.

In the present context however, we can say more.

\begin{proposition}\label{prop3.8}
We have $W\in M(A\ot B)$ and $\langle W, b\ot a\rangle =\langle
a,b\rangle$ for all $a\in A$ and $b\in B$. Also $W^{-1}\in M(A \ot B)$
and $\langle W^{-1},b \ot a\rangle=\langle S(a), b \rangle
=\langle a, S^{-1}(b) \rangle$.
\end{proposition}

These formulas need some explanation. The algebra $A\ot B$ acts on
$A\ot A$ in the obvious way (where we use the action $\la$ of $B$ and
$\pi$ of $A$ as defined in Definition~\ref{def3.1}). This action of $A\ot B$ is
non-degenerate and therefore extends to an action of the multiplier
algebra $M(A\ot B)$ in a unique way (see e.g.\ \cite{[D-VD-Z]}).
The first statement in the above
proposition says that the linear operator~$W$ is the action of some
multiplier in $M(A\ot B)$. For the second result, we observe that we
have a natural pairing of $A \ot B$ with $B\ot A$ (coming from the
pairing of $A$ with $B$) and that this pairing can be extended to a
bilinear map pairing $M(A\ot B)$ with $B\ot A$.

The fact that $W\in M(A\ot B)$ makes it possible to interpret certain
(well-known) formulas in another, nicer way. One can now show that
the formula
\[(\De \ot \iota)W=W_{13}W_{23}\]
makes sense in the algebra $M(A\ot A\ot B)$. For this we have to
observe that~$\mbox{$\De\ot\iota$}$ is a non-degenerate $^*$-homomorphism from
$A\ot B$ to $M(A\ot A)\ot B$, which is naturally imbedded in
$M(A\ot A\ot B)$, and hence has a unique extension to a unital
$^*$-homomorphism from $M(A\ot B)$ to $M(A\ot A\ot B)$.

The other equation, namely $\De(a)=W^{-1}(1\ot a)W$ is somewhat more
tricky. This equation is clear when considered as an equation of linear
maps from~$\mbox{$A\ot A$}$ to itself. But in these circumstances, it can also
be viewed as an equation in $M(A\ot C)$ where $C$ is the algebra
generated by $A$ and $B$, taking into account the commutation relations
of Lemma~\ref{lem3.2}. This algebra $C$ is called the {\it Heisenberg algebra}
(or sometimes the Heisenberg double) but,
as we mentioned before, its structure only depends on the
pairing of the linear spaces $A$ and $B$. Therefore we tend to think of
the Heisenberg algebra $C$ as an algebra with two special subalgebras
$A$ and $B$, sitting in the multiplier algebra $M(C)$, rather then of
the algebra itself.

Finally, there is the equation
\[W^{-1}=(S\ot\iota)W=(\iota\ot S^{-1})W.\]
This equation can be given a meaning using not only that
{$W\in M(A\ot B)$}~but also that $(a\ot 1)W(1\ot b)$ and $(1 \ot b)W(a\ot 1)$
are in fact elements in~$A\ot B$. If e.g.\ we apply $S\ot\iota$ to the
first of these two elements, we get formally,
because~$S$~is~an anti-homomorphism
$((S\ot\iota)W)(S(a)\ot b)$ and this is then~$W^{-1}(S(a)\ot b)$.

Using the fact that $W$ is (in a way) the duality, we can rewrite the
first formula of Proposition~\ref{prop3.4} as
\begin{align*} b&=F(a)=(\vfi\ot\iota)(W(a\ot 1)) \\
         a&=F^{-1}(b)=(\iota\ot\vfi)(W^{-1}(1\ot b)).
\end{align*}
Using these expressions, we have an easy way to obtain the Plancherel
formula. Indeed, using the same notations as above,
\begin{align*} \vfi(b^*b)&=\vfi(b^*((\vfi\ot\iota)W(a\ot 1))) \\
&=(\vfi\ot\vfi)((1\ot b^*)W(a\ot 1)) \\
&=\vfi((\iota\ot\vfi)((1\ot b^*)W)a)=\vfi(a^*a)
\end{align*}
We have used that $W$ is a unitary in the {$^*$-algebra} case and so
$W^*=W^{-1}$.

The Pentagon equation for $W$ (Lemma~\ref{lem3.7}) is in fact equivalent with
the Heisenberg commutation relations (Lemma~\ref{lem3.2}) given that $W$ is the
duality (cf.\ Proposition~\ref{prop3.8}). The relation with the Heisenberg
commutation rules becomes also more apparent in the following formula.

\begin{proposition}\label{prop3.9}
If we transform the map $W$ on $A\ot A$ to a map on $B\ot A$ using the
Fourier transform $F$ (i.e.\ $F_1$ in Definition~\ref{def3.3}) on the first
leg, we get the map
\[y\ot x\mapsto \langle S^{-1}(x_{(1)}),y_{(1)}\rangle y_{(2)} \ot
x_{(2)}.\]
\end{proposition}

Indeed, we had $W(x'\ot x)=S^{-1}(x_{(1)})x'\ot x_{(2)}$ and we know
from Proposition~\ref{prop3.5} that $F(ax')=\la(a)F(x')$ where
$\la(a)y=\langle a,y_{(1)}\rangle y_{(2)}$.

It should be observed that the inverse of the above map is
\[y\ot x\mapsto \langle x_{(1)},y_{(1)}\rangle y_{(2)} \ot
x_{(2)}\]
and that these two maps are indeed well-defined maps from $B\ot A$ to
itself (which is not completely obvious but follows from the previous
considerations, see also \cite{[D-VD]}).
Also remark that these two maps precisely govern the
Heisenberg commutation relations (see Lemma~\ref{lem3.2}).

We finish this section by another application of the above expressions.

\begin{proposition}
Let $A$ be finite-dimensional and denote by $\tr$ the trace on the algebra
$C$. Then, for some scalar $k\in \mathbb{C}$, we have
\[\vfi(a)=k\,\tr(\pi(a))\]
for all $a\in A$.
\end{proposition}

 \textbf{\emph{Proof.}}\,
We will show that the right hand side is left invariant and then we
get the formula by uniqueness of left invariant functionals. Observe
that $\tr(\pi(1))=\tr(1)\neq 0$.

So let $a\in B$. We will write $\sum u_i\ot v_i$ for $W$. We will
also drop $\pi$ and simply write $a$ for $\pi(a)$. Then we have
\begin{align*}
(\iota\ot\tr)\De(a)
&=(\iota\ot\tr)(W^{-1}(1\ot a)W) \\
&=(\iota\ot\tr)((S \ot \iota)W(1\ot a)W) \\
&=\sum(\iota\ot\tr)((1\ot a)(S(u_i)\ot 1)W(1\ot v_i))\\
&=\sum S((\iota\ot\tr)((1\ot a)W^{-1}(u_i\ot v_i)))\\
&=\tr(a)S(1)=\tr(a).
\end{align*}
\ebew\\

Observe that we have used the trace property and also the fact that in
this case $S^2=\iota$. By composing with the antipode, we see that
here the trace is also left invariant. The above argument can be
used to prove the existence of integrals on Hopf~{$^*$-algebras} with a
nice underlying {$^*$-algebra} structure.

In fact, the argument can be modified so that it also works
for any finite-dimensional Hopf algebra. In this case, one has to make
one modification because possibly $S^2\neq\iota$ and another one
because it might happen that the above functional is trivially $0$. We refer
to \cite{[K-VD2]}.

\section{Actions on Hilbert spaces}\label{sec4}

In the previous section, we essentially
did not use the $^*$-structure. In
this section, we will see the consequences of the positivity of the
integrals. We refer to \cite{[K-V2]}, \cite{[K-V3]} and \cite{[K-VD1]}
where more details can be found; see also \cite{[K-VD2]}.

A basic construction is the so-called {\it G.N.S.-construction}:

\begin{definition}
Let $\mathcal{H}$ denote the Hilbert space obtained by completing $A$ for
the scalar product given by $(x',x)\mapsto \vfi(x^*x')$. Denote by
$x\mapsto \eta(x)$ the canonical imbedding of $A$ in $\mathcal{H}$ so that
$\langle\eta(x'),\eta (x)\rangle=\vfi(x^*x')$ for all $x,x'\in A$.
Similarly, we will use $\hat{\mathcal{H}}$ and $\hat\eta$ for these objects
for the integral $\vfi$ on $B$.
\end{definition}

We now have the following results concerning the actions $\pi$ and
$\la$ as defined in the previous section.

\begin{proposition}
The maps $\eta(x)\mapsto \eta(ax)$ and $\eta(x)\mapsto
\langle S^{-1}(x_{(1)}),b\rangle \eta (x_{(2)})$ extend to bounded linear
maps on $\mathcal{H}$. We will use~$\pi(a)$ for the first one and~$\la(b)$
for the second one (in agreement with the notations used in the
previous section). We have $\pi(a^*)=\pi(a)^*$ and
$\lambda(b^*)=\lambda(b)^*$
whenever $a\in A$ and $b\in B$.
\end{proposition}

These last formulas are relatively easy to verify, using the left
invariance of $\vfi$. However, the boundedness of the maps $\pi(a)$ and
$\la(b)$ is not so obvious. We will give an argument later.

We have similar results for $B$. That is, we have bounded maps $\pi(b)$
and~$\la(a)$ on $\hat{\mathcal{H}}$ given by
$\pi(b)\hat\eta(y)=\hat\eta(by)$ and $\la(a)\hat\eta(y)=\langle a,y_{(1)}
\rangle \hat\eta (y_{(2)})$ whenever $a\in A$ and $b,y\in B$. Also here
$\pi(b^*)=\pi(b)^*$ and $\lambda(a^*)=\lambda(a)^*$.

Because of the `Plancherel formula' (cf.\ a remark after
Proposition~\ref{prop3.4}), we have the following:

\begin{proposition}
The Fourier transform $F$ is an isometry of $\mathcal{H}$ onto $\hat{\mathcal{H}}$
and it transforms the operators $\pi(a)$ and $\la(b)$ on $\mathcal{H}$ to
respectively $\la(a)$ and $\pi(b)$ on~$\hat{\mathcal{H}}$.
\end{proposition}

This clarifies in a way the boundedness of the operators $\lambda(a)$
and $\lambda(b)$, provided we know the boundedness of the operators
$\pi(a)$ and $\pi(b)$ already, but that is more obvious.\\

Now, we will see what we can say about the map $W$.

\begin{proposition}
The linear operator $W$ as defined from $A\ot A$ to itself, `extends'
to a unitary operator on $\mathcal{H}\ot \mathcal{H}$, still denoted by $W$.
\end{proposition}

So we have
\begin{align*} W(\eta(x)\ot \eta(x'))&=\eta(S^{-1}(x_{(1)}')x)\ot
\eta(x'_{(2)})\\
W^*(\eta(x)\ot \eta(x'))&=\eta(x_{(1)}'x)\ot
\eta(x'_{(2)}).
\end{align*} The unitary $W$ is called the {\it left regular
representation} of $(A,\De)$. If we take any other $x''\in A$, and if
we define a normal linear functional $\omega$ on $\mathcal{B}(\mathcal{H})$ (the
algebra of all bounded linear operators on $\mathcal{H}$) by $\omega(z)=
\langle z\eta(x'),\eta (x'')\rangle$, we see that $(\iota\ot \omega)
W^*=\vfi({x''}^*x'_{(2)})\pi(x'_{(1)})$. It is this type of formula that
can be used to argue that $\pi(a)$ is a bounded operator on $\mathcal{H}$.
And similarly for $\la(b)$.

Again, we can apply the Fourier transform $F$. If we only apply it on
the first leg, we get a unitary map $U$ from $\hat{\mathcal{H}}\ot \mathcal{H}$
to itself given by
\[U(\hat\eta(y)\ot \eta(x))=\langle S^{-1}(x_{(1)}), y_{(1)}\rangle
\,\, \hat\eta(y_{(2)})\ot \eta(x_{(2)}).\]
This is the Hilbert space version of the map given in Proposition~\ref{prop3.9}.

We can look at this formula in two ways. First, the right hand side is
\[\la(S^{-1}(x_{(1)}))\hat\eta(y)\ot \eta(x_{(2)}).\]
This is the form we get when we start from the formula with $W$ and
apply the Fourier transform $F$ on the first leg in the tensor
product.

We can also look at the right hand side as
\[\hat\eta(y_{(2)}) \ot \la(y_{(1)})\eta(x)\]
and when we apply the Fourier transform once more, now on the second
leg, we get the operator on $\hat{\mathcal{H}} \ot \hat{\mathcal{H}}$ given by
\[\hat\eta(y)\ot\hat\eta(y')\mapsto
\hat\eta(y_{(2)}) \ot \hat\eta (y_{(1)}y').\]
When we flip the two components in the tensor product
$\hat{\mathcal{H}} \ot \hat{\mathcal{H}}$, we precisely get the adjoint of the
unitary operator
$\hat W$ defined by
\[\hat W(\hat\eta(y)\ot\hat\eta(y'))= \hat\eta(S^{-1}(y'_{(1)})y')
\ot \hat\eta(y_{(2)}).\]
This is the left regular representation of $(B,\De)$. So we find the
that the Fourier transform carries $W$ into $\Sigma\hat W^*\Sigma$
where $\Sigma$ is used here to denote the flip
on~$\hat{\mathcal{H}}\ot \hat{\mathcal{H}}$.\\

We will come back to this operator later. Now, we want to apply the
Tomita-Takesaki theory. According to this theory, we get the following
(see e.g.\ \cite{[St]}). We denote by $M$ the weak operator closure of
$\pi(A)$ and by $M'$ the commutant of $M$.

\begin{proposition}
There exists a positive, self-adjoint, non-singular operator $\nabla$
and a unitary, conjugate linear involution $J$ on $\mathcal{H}$ such that
$\eta(x)\in \mathcal{D}(\nabla^{\frac12})$ and
$\eta(x^*)=J\nabla^{\frac12}\eta(x)$ for all $x\in A$. The space
$\eta(A)$ is a core for the domain of $\nabla^{\frac12}$.
We also have that $JMJ=M'$ and
$\nabla^{it}M\nabla^{-it}=M$ for all $t\in \mathbb{R}$.
\end{proposition}

The maps $\sigma_t$, defined by
$\sigma_t(x)=\nabla^{it}x\nabla^{-it}$ are called the {\it modular
automorphisms}. One can show that $\pi(a)$ for $a\in A$ is analytic for
this one-parameter group of automorphisms and that
$\sigma_{-i}(\pi(a))=\pi(\sigma(a))$ (see \cite{[K-VD1]}). Do not confuse
the automorphism $\sigma$ on $A$ with the one-parameter group of
automorphisms $(\sigma_t)$ on $M$.

Similarly, we have operators $\hat\nabla$ and $\hat J$ on $\hat{\mathcal{H}}$
satisfying $\hat\eta(y^*)=\hat J\hat\nabla^{\frac12}\hat\eta(y)$
whenever $y\in B$. They have the same properties w.r.t.\ the von
Neumann algebra $\hat M$, obtained by taking the weak closure of
$\pi(B)$ on $\hat{\mathcal{H}}$.\\

Let us now consider a very {\it special functional} on the Heisenberg
algebra $C$, introduced in the previous section. For a treatment of
this functional in the case of a general locally compact quantum group,
see \cite{[V-VD]}.

\begin{definition}\label{def4.6}
Define a linear functional~$f$ on the Heisenberg algebra~$C$ by
$f(ba)=\vfi(b)\vfi(a)$ for $a\in A$ and $b\in B$.
\end{definition}

\begin{lemma}\label{lem4.7}
This functional is positive and $f(a^*ba)=\vfi(b)\vfi(a^*a)$ for
all $a\in A$ {and~$b\in B$}.
\end{lemma}

This is easy to verify. Indeed
\[aba'=\langle a_{(1)},b_{(1)} \rangle b_{(2)}a_{(2)}a'\]
so that
\begin{align*}
f(aba')&=\langle a_{(1)},b_{(1)} \rangle \vfi(b_{(2)})\vfi(a_{(2)}a')\\
       &=\langle a_{(1)},1 \rangle \vfi(b)\vfi(a_{(2)}a')\\
       &=\vfi(b)\vfi(aa').
\end{align*}

Now, we want to apply the G.N.S.-construction for $f$ and relate the
Tomita-Takesaki data for $f$ with those for $\vfi$ on $A$ and $\vfi$
on $B$ as introduced above.

We first have the following.

\begin{proposition}
We can identify the Hilbert space $\mathcal{H}_f$ with $\hat{\mathcal{H}}\ot \mathcal{H}$ and
we get for the canonical imbedding
\[\eta_f(yx)=\hat\eta(y)\ot\eta(x)\]
whenever $x\in A$ and $y\in B$. The G.N.S.-representation $\pi_f$ is
given by
\begin{align*}
 \pi_f(a)&=\la(a_{(1)})\ot \pi (a_{(2)}) \\
 \pi_f(b)&=\pi(b)\ot 1
\end{align*}
when $a\in A$ and $b\in B$.
\end{proposition}

 \textbf{\emph{Proof.}}\,
The first statement is obvious as
\[f((yx)^*(yx))=f(x^*y^*yx)=\vfi(y^*y)\vfi(x^*x)\]
for all $x\in A$ and $y\in B$ by Lemma~\ref{lem4.7}. As
$\pi_f(b)\eta_f(yx)=\eta_f(byx)$, we see immediately that
$\pi_f(b)=\pi(b)\ot 1$ when $b\in B$. Finally, using the commutation
rules of Lemma~\ref{lem3.2}, we get
\begin{align*} \pi_f(a)\eta_f(yx)&=\eta_f(ayx)
         =  \langle a_{(1)}, y_{(1)} \rangle\,\, \eta_f(y_{(2)}a_{(2)}x)\\
    &= \langle a_{(1)}, y_{(1)} \rangle\,\, \hat\eta(y_{(2)}) \ot
          \pi(a_{(2)}) \eta(x) \\
    &=\la(a_{(1)})\hat\eta(y)\ot \pi (a_{(2)})\eta(x)
\end{align*}
whenever $a,x\in A$ and $b\in B$, giving the last formula.
\ebew\\

We know from the previous observations that the comultiplication $\De$
on $A$ is implemented by $W$ in the sense that $\Delta(a)=W^*(1\ot a)W$.
Consider this formula on $\mathcal{H} \ot \mathcal{H}$ and apply the Fourier
transform on the first leg. Then, we find
\[\pi_f(a)=U^*(1 \ot \pi(a))U\]
where $U$ is, as before, given by
\begin{align*}
U(\hat\eta(y)\ot \eta(x))&=\langle S^{-1}(x_{(1)}), y_{(1)} \rangle
\,\,\hat\eta (y_{(2)}) \ot \eta(x_{(2)}) \\
&=\hat\eta(y_{(2)})\ot \la(y_{(1)})\eta(x)\\
&=\la(S^{-1}(x_{(1)}))\hat\eta(y)\ot\eta (x_{(2)}).
\end{align*}

Next, we consider the polar decomposition $J_f\nabla_f^{\frac12}$ of
the operator $T_f$ defined as the closure of the map $\eta_f(z)\mapsto
\eta_f(z^*)$ with $z$ in the Heisenberg algebra $C$.

\begin{proposition}\label{prop4.9}
We have
\begin{align*} J_f &=(\hat J\ot J) U=U^* (\hat J \ot J) \\
         \nabla_f &=\hat\nabla \ot \nabla
\end{align*}
where $U$ and also $J$, $\hat J$, $\nabla$ and $\hat\nabla$ are as
above. Moreover, $U$ commutes with~$\nabla_f$.
\end{proposition}

 \textbf{\emph{Proof.}}\,
For $x\in A$ and $y\in B$ we have
\begin{align*} T_f\eta_f(yx)&=\eta_f(x^*y^*) \\
  &=\langle x_{(1)}^*, y_{(1)}^* \rangle \eta_f(y_{(2)}^*x_{(2)}^*)\\
  &=U^* \hat\eta(y^*) \ot \eta (x^*).
\end{align*}
It follows, by the uniqueness of the polar decomposition, that
\begin{align*} J_f &=U^* (\hat J \ot J) \\
         \nabla_f &=\hat\nabla \ot \nabla.
\end{align*}
Using the properties of the operators involved, we find easily that
also $ J_f =(\hat J\ot J) U$
and that $U$ commutes with $\nabla_f$.
\ebew\\

This result for general locally compact quantum groups can be found
e.g.\ in \cite{[K-V3]}, Corollary~2.2. And as in that same paper, some
nice conclusions can be drawn from the previous result.

For this we need to introduce some notations. We have used $M$ and
$\hat M$ to denote the von Neumann algebras on $\mathcal{H}$ and
$\hat{\mathcal{H}}$ generated by $\pi(a)$ with $a\in A$ and $\pi(b)$
with $b\in B$ respectively. We will use $\hat N$ and $N$ to denote the
von Neumann algebras on $\mathcal{H}$ and $\hat{\mathcal{H}}$
generated by $\lambda(b)$ with $b\in B$ and $\lambda(a)$ with $a\in A$
respectively.

Observe that the pair $(M,\hat N)$ acts on~$\mathcal{H}$ and that the pair
$(N, \hat M)$ acts on~$\hat{\mathcal{H}}$. The Fourier transform  sends $M$
to $N$ and $\hat N$ tot $\hat M$ (because of Proposition~\ref{prop3.5}).

The operator $U$ has its first leg in $N$ and its second leg in $\hat
N$ (see e.g.\ the formulas for $U$ given before Proposition~\ref{prop4.9}).
Because the two legs of $U$ generate these von Neumann algebras, we
get, as a consequence of Proposition~\ref{prop4.9}, the following (see e.g.\
Proposition~2.1 in \cite{[K-V3]}).

\begin{proposition}
\begin{align*}
\hat\nabla^{it} N \hat\nabla^{-it} &=N \qquad\qquad\qquad
             \hat J N \hat J = N & \\
             \nabla^{it}\hat N \nabla^{-it} &=\hat N \qquad\qquad\qquad
             J \hat N J = \hat N. &
\end{align*}
\end{proposition}

We can be more specific. Define linear maps $R$ and $\hat R$ on $N$ and
$\hat N$ respectively by $R(x)=\hat J x^* \hat J$ and $\hat R(y)=J
y^*J$. Also define one-parameter groups of automorphisms $(\tau_t)$ and
$\hat \tau_t$ on $N$ and $\hat N$ respectively by
$\tau_t(x) = \hat\nabla^{it} x \hat\nabla^{-it}$ and
$\hat\tau_t(y) = \nabla^{it} y \nabla^{-it}$. Then we have the
following.

\begin{proposition}
We have $\lambda(a)\in\mathcal{D}(\tau_{-\frac{i}{2}})$ and
$\lambda(S(a))=R\tau_{-\frac{i}{2}}(\lambda(a))$ for all $a\in A$.
Similarly, $\lambda(b)\in\mathcal{D}(\hat\tau_{-\frac{i}{2}})$ and
$\lambda(S(b))=\hat R\hat\tau_{-\frac{i}{2}}(\lambda(b))$ for all
$b\in B$.
\end{proposition}

We will not give a rigorous proof, just indicating why these formulas
are true.

Take first $a\in A$ and $y\in B$. Then
\begin{align*} \hat J\hat\nabla^{\frac12} \lambda(a)\hat\eta(y)
&=\hat J\hat\nabla^{\frac12}\langle a,y_{(1)} \rangle \hat\eta (y_{(2)}) \\
&=\langle a,y_{(1)} \rangle^- \hat\eta (y_{(2)}^*) \\
&=\langle S(a)^*,y_{(1)}^* \rangle \hat\eta (y_{(2)}^*)\\
&=\lambda(S(a)^*) \hat\eta (y^*)\\
&=\lambda(S(a)^*) \hat J\hat\nabla^{\frac12}  \hat\eta (y).
\end{align*}
We get formally
\[\lambda(S(a)^*)= \hat J\hat\nabla^{\frac12} \lambda(a)
\hat\nabla^{-\frac12}\hat J\]
and so $\lambda(S(a))=R\tau_{-\frac{i}{2}}(\lambda(a)).$

Similarly,
take $b\in B$ and $x\in A$. Then
\begin{align*} J\nabla^{\frac12} \lambda(b)\eta(x)
&= J\nabla^{\frac12}\langle S^{-1}(x_{(1)}),b \rangle \eta (x_{(2)}) \\
&=\langle x_{(1)},S(b) \rangle^- \eta (x_{(2)}^*) \\
&=\langle x_{(1)}^*,b^* \rangle \eta (x_{(2)}^*)\\
&=\langle S^{-1}(x_{(1)}^*),S^{-1}(b^*) \rangle \eta (x_{(2)}^*)\\
&=\lambda(S^{-1}(b^*) \eta (x^*)\\
&=\lambda(S(b)^*) J\nabla^{\frac12}  \eta (x).
\end{align*}
Again formally
\[\lambda(S(b)^*)= J\nabla^{\frac12} \lambda(b)
\nabla^{-\frac12} J\]
and so $\lambda(S(b))=\hat R\hat\tau_{-\frac{i}{2}}(\lambda(b)).$\\

As we see from the above, the use of the functional $f$, as defined in
Definition~\ref{def4.6} (and studied in more detail in \cite{[V-VD]}), in
combination with the use of the Fourier transform, admits a certain way
of understanding the very basic formulas obtained by Kustermans and Vaes
in \cite{[K-V2]} and \cite{[K-V3]}. The remaining step to arrive really
to their framework is to identify the spaces $\mathcal{H}$ and
$\hat{\mathcal{H}}$ so that the Fourier transform $F$ becomes the
identity map. Then the pairs $(M,\hat N)$ and $(N,\hat M)$ will become
the same (namely $(M,\hat M)$ in their terminology). We believe that it
has some advantage to wait for this identification until some level of
understanding has been obtained (as we did in this paper).

\end{document}